\documentclass[12pt]{article}
\usepackage{amsmath, amsfonts, amssymb, amsthm} 
\usepackage{enumerate}
\newtheorem{thm}{Theorem}[section]
\newtheorem{cor}[thm]{Corollary}
\newtheorem{lemma}[thm]{Lemma}
\newtheorem{prop}[thm]{Proposition}

\newtheorem{remark}[thm]{Remark}
\newtheorem{notation}[thm]{Notation}
\newtheorem{example}[thm]{Example}
\newtheorem{remarks}[thm]{Remarks}
\newtheorem{definition}[thm]{Definition}

\theoremstyle{definition}
\newtheorem{defin}[thm]{Definition}
\newtheorem{rem}[thm]{Remark}
\newtheorem{exa}[thm]{Example}

\numberwithin{equation}{section}

\def\R{\mathbb R} \def\Z{\mathbb Z} \def\C{\mathbb C} 

\def\N{\mathbb N}
\def\C{{\mathbb C}} 
\def\Q{\mathbb Q}

\def\<{\,<\!}
\def\>{\!>\,}
\def\Re{{\rm Re\,}}

\usepackage{graphics}
\usepackage{epsf}
\usepackage{epsfig}

\begin{document}

\title{ Convergents as approximants in continued fraction expansions of complex numbers with Eisenstein integers} 

\author{S.G. Dani}
\date{}
\maketitle


\renewcommand{\thefootnote}{}



\renewcommand{\thefootnote}{\arabic{footnote}}
\setcounter{footnote}{0}

\section{Introduction}

It is well-known that the convergents defined in terms of the simple continued fraction expansion of a 
real number are ``best approximants", in the sense that if $t\in \R$ and  $\{p_n/q_n\}$ is the 
corresponding sequence of convergents  then for any $n\in \N$ and  $1\leq q \leq q_n$, we have  $|qt-p|\geq |q_nt-{p_n}|$ for all $p\in \Z$ (see \cite{HW},Theorem~181, or \cite{NZ}, 
Theorem~7.13, for instance). 
In this paper we shall be concerned with the analogous issue  of  comparing  
 $|qz-p|$ with $|q_nz_n-p_n|$ (see below),   for continued fraction expansions of complex numbers. In this respect 
Hensley \cite{Hens} considered the continued fractions expansions of $z\in \C$ in terms of 
 the ring $\frak G$ of Gaussian integers,  defined via the nearest integer algorithm, and showed that if $\{p_n/q_n\}$ is the corresponding sequence of convergents   of  $z\in \C$, then  for any $n\in \N$ and $p,q\in \frak G$, such that  $1\leq  |q|\leq |q_n|$, we have $|qz-p|\geq \frac 15|q_nz - {p_n}|$. Here we prove an analogous result for the ring of Eisenstein integers, namely $\Z[\omega]$, where $\omega$ is a nontrivial cube root of unity. 

\begin{thm}\label{new-main}
Let $\frak E$ be  the ring of Eisenstein integers.  Let 
 $z\in \C$ and 
$\{a_n\}$  be a continued fraction expansion of $z$ over $\frak E$ with respect to the nearest integer algorithm and   let $\{p_n/q_n\}$, 
be the corresponding sequence of convergents.  
Then for any $q\in \frak E$ such that $1\leq |q| \leq |q_n|$ and  any $p\in \frak E$, $$|qz-{p}|\geq \frac 12 |q_nz-p_n|.$$   
\end{thm}

As an application of Theorem~\ref{new-main}, we obtain the following result analogous to the classical characterisation of badly approximable numbers. 

\begin{definition} {\rm  We say that $z\in \C$ is {\it badly approximable} with respect to $\frak E$ if there exists $\delta >0$ such that for all $p,q\in \frak E$, 
$q\neq 0$, $|z-\frac pq|\geq \delta/|q|^2$. }
\end{definition}

\begin{cor}\label{bad}
Let $\frak E$ be the ring of Eisenstein integers and  $K$ be 
the quotient field of $\frak E$. Let $z\in \C\backslash K$ and  $\{a_n\}_{n=0}^\infty$ be the continued fraction expansion of $z$ with respect to the nearest integer algorithm on $\frak E$. Then $z$ is badly approximable with respect to $\frak E$ if and only if $\{|a_n|\}_{n=0}^\infty$ is bounded. 

\end{cor} 

\section{Prelininaries}

Let $\frak E$ be the ring of Eisenstein integers. Then $\frak E$ can also be realised as $\Z[\rho]$, where $\rho = \frac 12 +\frac i2\sqrt 3$, which is a $6$th root of unity. 
We recall that the nearest integer algorithm over $\frak E$ is by defined a map $f: \C\to \frak E$ such that for all $z\in \C$, $|z-f(z)|\leq |z-a|$ for all $a\in \frak E$; the condition determines $f(z)$ uniquely for $z$ in the complement of a countable set of lines, while for other points, which are equidistant from distinct elements of $\frak E$, there can be multiple choices; and we shall call any $f$ as above a nearest integer algorithm (we may nevertheless refer to it as ``the nearest integer algorithm", as the multiple choices, for points for which  they are available, do not play any role in our results.). 

For any $z\in \C$  we get two sequences $\{a_n\}_{n=0}^m$
 and $\{z_n\}_{n=0}^m$, where $m$ is either a nonnegative integer or $\infty$, as follows: 
we set $z_0=z$,  and having defined $z_0, \dots, z_n$ for some $n\geq 0$ we set $z_{n+1}=(z_n-f(z_n))^{-1}$ if $f(z_n)\neq z_n$, and terminate the sequence, choosing $m=n$, if $f(z_n)=z_n$, and define $a_n= f(z_n)$ for all $n=0,\dots , m$. The sequence $\{a_n\}_{n=0}^m$ is called the {\it continued fraction expansion of $z$}, $a_n$ are called the {\it partial quotients} of the expansion, and $\{z_n\}_{n=0}^m$ is called the corresponding {\it iteration sequence}. 

With the continued fraction expansion  $\{a_n\}_{n=0}^m$, of $z\in \C$, we associate  two
sequences $\{p_n\}_{n=-1}^\infty$ and $\{q_n\}_{n=-1}^\infty$  defined 
recursively by the relations
$$ p_{-1}=1, p_0=a_0, p_{n+1}=a_{n+1}p_n+p_{n-1}, \; \mbox{for all }
n\geq 0, \mbox{ and}
$$
$$
 q_{-1}=0, q_0=1, q_{n+1}=a_{n+1}q_n+q_{n-1}, \; \mbox{
  for all } n \geq 0.
$$
The pair of sequences $\{p_n\}_{n=-1}^\infty$, $\{q_n\}_{n=-1}^\infty$ is called the $\cal Q$-pair corresponding to the expansion. It is known that $q_n\neq 0$ for all $n\geq 1$ and if $m$ is infinite
$p_n/q_n\to z$ as $n\to \infty$; see \cite{D-Eisen}, for instance; $p_n/q_n$ are called the {\it convergents} of the expansion. We recall here the following result from \cite{D-Eisen}.

\begin{thm}\label{ratios}
Let the notation be as above. Then for all $n\geq 1$, $|q_{n+1}/q_{n-1}| \geq \frac 32.$
\end{thm}

Apart from the notation as above we shall also set, for all $n\geq 1$, $r_n=q_{n-1}/q_n$.  We 
note that Theorem~\ref{ratios} has the following obvious consequence. 

\begin{cor}\label{cor:ratios}
Let the notation be as above. Then for all $n\geq 1$, either $|r_n|\leq \sqrt \frac 23$, or $|r_{n+1}|\leq \sqrt \frac 23$.
\end{cor}

For the proof of Theorem~\ref{new-main} we recall also the following standard properties of the sequences associated with the  continued fraction expansions. 

\begin{prop}\label{prelim} For all $n\geq 0$  the following statements hold : i) $|z_{n+1}|\geq \sqrt 3$, 
ii) $|q_n|<|q_{n+1}|$,  
$$ iii) \ \ z=\displaystyle{\frac {p_nz_{n+1}+p_{n-1}}{q_nz_{n+1}+q_{n-1}}},\ \   
iv) \ \ \displaystyle |{q_n}z-{p_n}| =|\frac 1{(z_{n+1}q_n+q_{n-1})}|.$$
\end{prop}

\proof Assertion (i)  follows from the fact that $z_{n+1}^{-1}=(z_n-a_n)$ is contained in the ball or radius $1/\sqrt 3$, in fact in the hexagon with vertices at $\rho^ki/\sqrt 3$, $k=0, \dots, 5$ (see \cite{D-Eisen}, \S\,5. For Assertion~(ii) we refer the reader to \cite{L}, or to \cite{D-Eisen}, Theorem~5.1, where it is proved in a general setting. Assertion~(iii) and~(iv) are general facts about continued fraction expansions, that follow from straightforward manipulations on the recurrence relations; see \cite{D-Eisen}, Proposition~2.1(iii), and the proof of Proposition~2.1(iv),  for ready reference). \qed

\section{On the quality of convergents as approximants}

Through the section we fix $z\in \C$ and follow the notation as in the last section associated with the continued fraction expansion of $z$ with respect to a nearest integer algorithm. We shall  now compare $|qz-p|$ with $|q_nz-{p_n}|$ for $q$ such that $1\leq |q|\leq |q_n|$,  and prove Theorem~\ref{new-main}, and then Corollary~\ref{bad}. 

\noindent{\it Proof of Theorem~\ref{new-main}}:
Let $n\in \N$ and $p,q\in \frak E$  be given such that $1\leq |q|\leq |q_n|$. 
Since $p_nq_{n-1}-q_np_{n-1}=(-1)^{n-1}$, it follows that 
there
exist $\alpha, \beta \in \frak E$ such that $\left (\begin{matrix} p \\ q\end{matrix}\right) =\alpha \left (\begin{matrix} p_{n} \\ q_{n} \end{matrix}\right)+\beta \left (\begin{matrix} p_{n-1} \\ q_{n-1} \end{matrix}\right)$; thus 
$p=\alpha p_{n} +\beta {p_{n-1}}$ and $q=\alpha q_{n} +\beta {q_{n-1}}$. 
Now, by Proposition~\ref{prelim} we have $z=\displaystyle{\frac {z_{n+1}p_n+p_{n-1}}{z_{n+1}q_n+q_{n-1}}}$, and hence
$$|qz-p|=|q||z-\frac pq| =|q||\frac {z_{n+1}p_n+p_{n-1}}{z_{n+1}q_n+q_{n-1}}-\frac {\alpha p_n+\beta p_{n-1}}{\alpha q_n+\beta q_{n-1}}|= |\frac {\beta z_{n+1} -\alpha}{z_{n+1}q_n+q_{n-1}}|, $$ using that $|p_nq_{n-1}-q_np_{n-1}|=1$ and $\alpha q_n+\beta q_{n-1}=q$. Also, by Proposition~\ref{prelim}(iv)
$$|{q_n}z-{p_n}| =|\frac 1{z_{n+1}q_n+q_{n-1}}|,$$
and hence $$|qz-p|= |\beta z_{n+1} -\alpha||{q_n}z-{p_n}|.$$
To prove the theorem it therefore suffices to show that $|\beta z_{n+1} -\alpha|\geq \kappa $ for all $\kappa <\frac 12$.  Let $\kappa <\frac 12$ be given. We 
shall suppose that $|\beta z_{n+1} -\alpha| <  \kappa$ and arrive at a contradiction. 

First suppose that, if possible, $\beta =0$. Then by the assumption as above $|\alpha |<\kappa <1$ and since
$\alpha \in \frak E$ we get  that $\alpha =0$. But in turn this implies that $q=\alpha q_n+\beta q_{n+1}=0$, contrary to the hypothesis. Hence $\beta \neq 0$.  Now suppose  that $|\beta| =1$. Then replacing $p$ and $q$ by their multiples by a fixed unit in $\frak E$ we may assume $\beta =1$. Then we have $|z_{n+1} -\alpha|< \kappa$. Since $\kappa<\frac 12$ this  implies that
$\alpha =a_{n+1}$ and thus $q=\alpha q_n+\beta q_{n-1} =a_{n+1} q_n+q_{n-1}=q_{n+1}$, which is a contradiction since by hypothesis 
$|q|\leq |q_n| <|q_{n+1}|$, where the last inequality is as in Proposition~\ref{prelim}(ii). Thus $|\beta |>1$, and since $\beta \in \frak E$ we get that  $|\beta|\geq \sqrt 3$. 

Since $ |\alpha q_n+\beta q_{n-1}| = |q|\leq |q_{n}|$, dividing by $|q_n|$ we get 
$ |\alpha +\beta r_n|\leq 1$. Since by assumption $|\beta z_{n+1}-\alpha|< \kappa$, we get that  $$ |\beta ||z_{n+1}+r_n|=|\beta z_{n+1}+\beta r_n|\leq |\beta z_{n+1}-\alpha|+ |\alpha +\beta r_n|<\kappa +1.$$ Recall that by Corollary~\ref{cor:ratios} either 
$|r_n|$ or $|r_{n-1}|$ is at most $\sqrt {\frac 23}$. Suppose $|r_n|\leq \sqrt {\frac 23} $. Then, recalling that $|\beta |\geq \sqrt 3$ and $|z_{n+1}|\geq \sqrt 3 $ (cf. Proposition~\ref{prelim}(i)),  
we get $$\kappa +1> |\beta ||z_{n+1}+r_n|\geq  |\beta |(|z_{n+1}|-|r_n|)\geq \sqrt 3 \left (\sqrt 3 - \sqrt {\frac 23} \right )=3-\sqrt 2 >\frac 32,$$ which however contradicts the choice of $\kappa$.  

Now suppose that $|r_n|>\sqrt {\frac 23}$ and $|r_{n-1}|\leq \sqrt {\frac 23}$. We have $z_{n+1}=1/(z_n-a_n)$ and $r_n=
1/(a_n+r_{n-1})$, and hence $$|z_{n+1}+r_n|=|\frac 1{z_n-a_n}+\frac 1{a_n+r_{n-1}}|=\frac {|z_n+r_{n-1}|}{|(z_n-a_n)(a_n+r_{n-1})|}=|z_n+r_{n-1}||z_{n+1}||r_n|.$$ 
As $|\beta |\geq \sqrt 3$, $|z_n|\geq \sqrt 3$, $|r_{n-1}|\leq \sqrt {\frac 23}$, $|z_{n+1}|\geq \sqrt 3$ and $|r_n|> \sqrt {\frac 23}$, it follows that  $$|\beta ||z_{n+1}+r_n|=|\beta| |z_n+r_{n-1}||z_{n+1}||r_n|\geq \sqrt 3\left (\sqrt 3-\sqrt {\frac 23}\right) \sqrt 3 \sqrt {\frac 23} > \frac 32.$$ Thus  we get $\kappa +1>  \frac 32$, 
again contradicting the choice of $\kappa$. Therefore $|\beta z_{n+1}-\alpha|\geq \kappa$, as sought to be proved. \qed

\bigskip
\noindent {\it Proof of Corollary~\ref{bad}}:   Since $|z_{n+1}-a_{n+1}|<1$ and $|q_{n-1}|<|q_n|$, we have  $$|a_{n+1}|-2 \leq |z_{n+1}+\frac {q_{n-1}}{q_n}|\leq |a_{n+1}|+2,$$ for all $n$. 
Since, by Proposition~\ref{prelim}(iv), we have 
 $|q_n||q_nz-{p_n}| = (|z_{n+1}+\frac {q_{n-1}}{q_n}|)^{-1}$, it follows that 
if  $\{|a_n|\}$ is unbounded then for any $\delta >0$ there exists $n$ such that $|q_n||q_nz-{p_n}|< \delta $, so $z$ is not badly approximable. 

Now suppose that $\{|a_n|\}$ 
 is bounded, say $|a_n|\leq M$ for all $n$. Then $$|q_nz-{p_n}|=|q_n|^{-1}(|z_{n+1}+\frac {q_{n-1}}{q_n}|)^{-1}\geq \frac 1{(M+2)|q_n|}$$ for all $n$. Now  let $p,q\in \frak E$ with $q\neq 0$ be arbitrary, and let $n\in \N$ be such that 
 $|q_{n-1}|\leq |q|\leq |q_n|$.  Since $q_n=a_nq_{n-1}+q_{n-2}$ and $|q_{n-2}|<|q_{n-1}|$, 
 we then have $|q_n|\leq (M+1)|q_{n-1}| \leq (M+1)|q|$. By Theorem~\ref{new-main} we have $|qz-p|\geq \frac 12 |{q_n}z-{p_n}|$ and hence $$|q||qz- p|\geq \frac 12|q|{|q_nz- p_n|}\geq \frac 12\frac {|q|}{(M+2)|q_n|}\geq \frac1{2(M+2)^2}, $$ and hence $|z-\frac pq|\geq \delta |q|^{-2}$, with $\delta =\frac 12(M+2)^{-2}$. Since this holds for all $p,q\in \Gamma$, $q\neq 0$ we get that $z$ is badly approximable. \qed

\begin{remark}
{\rm An argument as in the proof of Theorem~\ref{new-main} does not yield a similar result  in the case of the ring $\frak G$ of Gaussian integers (for which the corresponding result is proved in \cite{Hens} with the constant $\frac 15$ in place of $\frac 12$), since in getting a lower bound for $\kappa +1$ as in the proof we would only have at our disposal the estimates $|z_{n+1}|\geq \sqrt 2$, and the resulting bound may be seen to  be not good enough. }
 
\end{remark}

\vskip5mm
\begin{flushleft}
S.\,G.\,Dani\\
 Department of Mathematics\\
 Indian Institute of Technology Bombay \\
 Powai, Mumbai 400076\\ 
 India\\
 
 \smallskip
 E-mail: {\tt sdani@math.iitb.ac.in}
 \end{flushleft}

\end{document}